\documentclass[12pt]{amsart}
\usepackage{amsmath}
\usepackage{mathtools}
\usepackage{amssymb}
\usepackage{graphicx}
\usepackage{xcolor}
\usepackage{tikz}
\usepackage{tkz-graph}
\usetikzlibrary{cd}
\usetikzlibrary{knots}
\usetikzlibrary{calc}
\usetikzlibrary{decorations,decorations.pathreplacing,decorations.markings}
\tikzset{string/.style={ultra thick}}
\tikzset{smallstring/.style={thick,scale=0.75,every node/.style={transform shape}}}
\tikzstyle{every picture}+=[baseline={([yshift=-.5ex]current bounding box.center)}]

\tikzstyle{vertex}=[circle, draw, inner sep=0pt, minimum size=6pt]
\newcommand{\vertex}{\node[vertex]}
\usepackage[all]{xy}
\graphicspath{ {graphs/} }
\usetikzlibrary{decorations.text}

\usepackage{graphicx}
\usepackage{fullpage}
\usepackage[T1]{fontenc}
\usepackage{microtype}
\usepackage{ifthen}
\usepackage{libertine}
\usepackage[libertine]{newtxmath}
\usepackage[all]{xy}
\usepackage{pdflscape}
\usepackage{array}
\usepackage{multirow}
\usepackage{rotating}
\usepackage{mathtools}

\usepackage{xcolor}
\definecolor{dark-red}{rgb}{0.7,0.25,0.25}
\definecolor{dark-blue}{rgb}{0.15,0.15,0.55}
\definecolor{medium-blue}{rgb}{0,0,0.65}
\definecolor{DarkGreen}{RGB}{0,150,0}

\usepackage[pdftex,plainpages=false,hypertexnames=false,pdfpagelabels,breaklinks]{hyperref}
\hypersetup{
   colorlinks, linkcolor={purple},
   citecolor={medium-blue}, urlcolor={medium-blue}
}

\usepackage[backend=bibtex,style=alphabetic,doi=false,isbn=false,url=false,minnames=6,maxnames=6]{biblatex}
\renewbibmacro{in:}{%
  \ifentrytype{article}{}{\printtext{\bibstring{in}\intitlepunct}}}
\renewbibmacro*{volume+number+eid}{%
  \printtext{vol.}
  \printfield{volume}%
  \setunit*{\addnbspace}% NEW (optional); there's also \addnbthinspace
  \printfield{number}%
  \setunit{\addcomma\space}%
  \printfield{eid}}
\DeclareFieldFormat[article]{number}{\mkbibparens{#1}}
\usepackage{silence}
\allowdisplaybreaks
\newcommand{\doi}[1]{\href{https://dx.doi.org/#1}{{\small DOI:#1}}}

\newcommand{\googlebooks}[1]{(preview at \href{https://books.google.com/books?id=#1}{google books})}

\newcommand{\numdam}[1]{}

\theoremstyle{plain}
\numberwithin{equation}{section}

\theoremstyle{remark}

\theoremstyle{definition}
\newcommand{\sslash}{\mathbin{/\mkern-6mu/}}

\tikzstyle{mid>}=[decoration={markings, mark=at position 0.5 with {\arrow{>}}}, postaction={decorate}]
\tikzstyle{mid<}=[decoration={markings, mark=at position 0.5 with {\arrow{<}}}, postaction={decorate}]
\tikzstyle{upper>}=[decoration={markings, mark=at position 0.8 with {\arrow{>}}}, postaction={decorate}]
\tikzstyle{upper<}=[decoration={markings, mark=at position 0.8 with {\arrow{<}}}, postaction={decorate}]
\tikzstyle{lower>}=[decoration={markings, mark=at position 0.25 with {\arrow{>}}}, postaction={decorate}]
\tikzstyle{lower<}=[decoration={markings, mark=at position 0.25 with {\arrow{<}}}, postaction={decorate}]
\tikzstyle{tag}=[decoration={markings, mark=at position 0.5 with {\arrow{Rays[n=1]}}}, postaction={decorate}]

\newcommand{\iso}{\cong}

\DeclareMathOperator{\Hom}{Hom}
\DeclareMathOperator{\Inv}{Inv}

\newcommand{\sVec}{{\mathsf {sVec}}}
\newcommand{\Sem}{{\mathsf {Sem}}}
\newcommand{\Rep}{{\mathsf {Rep}}}
\renewcommand{\Vec}{{\mathsf {Vec}}}

\newcommand{\ad}{{\operatorname{Ad}}}

\def\semicolon{;}
\def\applytolist#1{
    \expandafter\def\csname multi#1\endcsname##1{
        \def\multiack{##1}\ifx\multiack\semicolon
            \def\next{\relax}
        \else
            \csname #1\endcsname{##1}
            \def\next{\csname multi#1\endcsname}
        \fi
        \next}
    \csname multi#1\endcsname}

\def\calc#1{\expandafter\def\csname c#1\endcsname{{\mathcal #1}}}
\applytolist{calc}QWERTYUIOPLKJHGFDSAZXCVBNM;
\def\bbc#1{\expandafter\def\csname bb#1\endcsname{{\mathbb #1}}}
\applytolist{bbc}QWERTYUIOPLKJHGFDSAZXCVBNM;
\def\bfc#1{\expandafter\def\csname bf#1\endcsname{{\mathbf #1}}}
\applytolist{bfc}QWERTYUIOPLKJHGFDSAZXCVBNM;

\makeatletter
\newlength{\L@UnitsRaiseDisplaystyle}
\newlength{\L@UnitsRaiseTextstyle}
\newlength{\L@UnitsRaiseScriptstyle}
\DeclareRobustCommand*{\@UnitsNiceFrac}[3][]{%
  \ifthenelse{\boolean{mmode}}{%
    \settoheight{\L@UnitsRaiseDisplaystyle}{%
      \ensuremath{\displaystyle#1{M}}%
    }%
    \settoheight{\L@UnitsRaiseTextstyle}{%
      \ensuremath{\textstyle#1{M}}%
    }%
    \settoheight{\L@UnitsRaiseScriptstyle}{%
      \ensuremath{\scriptstyle#1{M}}%
    }%
    \settoheight{\@tempdima}{%
      \ensuremath{\scriptscriptstyle#1{M}}%
    }%
    \addtolength{\L@UnitsRaiseDisplaystyle}{%
      -\L@UnitsRaiseScriptstyle%
    }%
    \addtolength{\L@UnitsRaiseTextstyle}{%
      -\L@UnitsRaiseScriptstyle%
    }%
    \addtolength{\L@UnitsRaiseScriptstyle}{-\@tempdima}%
    \mathchoice
      {%
        \raisebox{\L@UnitsRaiseDisplaystyle}{%
          \ensuremath{\scriptstyle#1{#2}}%
        }%
      }%
      {%
        \raisebox{\L@UnitsRaiseTextstyle}{%
          \ensuremath{\scriptstyle#1{#2}}%
        }%
      }%
      {%
        \raisebox{\L@UnitsRaiseScriptstyle}{%
          \ensuremath{\scriptscriptstyle#1{#2}}%
        }%
      }%
      {%
        \raisebox{\L@UnitsRaiseScriptstyle}{%
          \ensuremath{\scriptscriptstyle#1{#2}}%
        }%
      }%
    \mkern-2mu{\sslash}\mkern-1mu%
    \bgroup
      \mathchoice
        {\scriptstyle}%
        {\scriptstyle}%
        {\scriptscriptstyle}%
        {\scriptscriptstyle}%
      #1{#3}%
    \egroup
  }%
  {%
    \settoheight{\L@UnitsRaiseTextstyle}{#1{M}}%
    \settoheight{\@tempdima}{%
      \ensuremath{%
        \mbox{\fontsize\sf@size\z@\selectfont#1{M}}%
      }%
    }%
    \addtolength{\L@UnitsRaiseTextstyle}{-\@tempdima}%
    \raisebox{\L@UnitsRaiseTextstyle}{%
      \ensuremath{%
        \mbox{\fontsize\sf@size\z@\selectfont#1{#2}}%
      }%
    }%
    \ensuremath{\mkern-2mu}{\sslash}\ensuremath{\mkern-1mu}%
    \ensuremath{%
      \mbox{\fontsize\sf@size\z@\selectfont#1{#3}}%
    }%
  }%
}

\newcommand{\ncup}{\begin{tikzpicture}[baseline={([yshift=-.5ex]current bounding box.center)}]
\begin{scope}[very thick, every node/.style={sloped,allow upside down}]
\draw [thick] (0,0) arc [radius=-.25,start angle = 180, end angle = 0];
\end{scope}
\end{tikzpicture}}
\newcommand{\ncap}{\begin{tikzpicture}[baseline={([yshift=-.5ex]current bounding box.center)}]
\begin{scope}[very thick, every node/.style={sloped,allow upside down}]
\draw [thick] (0,0) arc [radius=.25,start angle = 180, end angle = 0];
\end{scope}
\end{tikzpicture}}

\newcommand{\rcup}{\begin{tikzpicture}[baseline={([yshift=-.5ex]current bounding box.center)}]
\begin{scope}[very thick, every node/.style={sloped,allow upside down}]
\draw [thick,mid>] (0,0) arc [radius=-.25,start angle = 180, end angle = 0];
\end{scope}
\end{tikzpicture}}
\newcommand{\lcup}{\begin{tikzpicture}[baseline={([yshift=-.5ex]current bounding box.center)}]
\begin{scope}[very thick, every node/.style={sloped,allow upside down}]
\draw [thick,mid<] (0,0) arc [radius=-.25,start angle = 180, end angle = 0];
\end{scope}
\end{tikzpicture}}

\newcommand{\lcap}{\begin{tikzpicture}[baseline={([yshift=-.5ex]current bounding box.center)}]
\begin{scope}[very thick, every node/.style={sloped,allow upside down}]
\draw [thick,mid>] (0,0) arc [radius=.25,start angle = 180, end angle = 0];
\end{scope}
\end{tikzpicture}}
\newcommand{\rcap}{\begin{tikzpicture}[baseline={([yshift=-.5ex]current bounding box.center)}]
\begin{scope}[very thick, every node/.style={sloped,allow upside down}]
\draw [thick,mid<] (0,0) arc [radius=.25,start angle = 180, end angle = 0];
\end{scope}
\end{tikzpicture}}

\newcommand{\JW}[1]{f^{(#1)}}

\newcommand{\xyJonesWenzlIdempotentSixPlusOne}[1]{\xybox{(0,-9.)*\xybox{*\xybox{
(1.5,0);(19.5,18) **i@{-},(3,18);(3,14.4) **@{-}, (3,3.6);(3,0) **@{-},(6,18);(6,14.4) **@{-}, (6,3.6);(6,0) **@{-},(9,18);(9,14.4) **@{-}, (9,3.6);(9,0) **@{-},(15,18);(15,14.4) **@{-}, (15,3.6);(15,0) **@{-},(18,18);(18,14.4) **@{-}, (18,3.6);(18,0) **@{-},(12,16.2) *h+++{\cdots},(12,1.8) *+++{\cdots},(2.25,3.6);(18.75,14.4) **\frm{-}, (10.5,9.) *+++{#1}}!R!(-1.5,0)*\xybox{
(1.5,0);(4.5,18) **i@{-},(3,18);(3,0)  **\crv{(3,5.4)&(3,12.6)},
}!R!(-1.5,0)}}}\newcommand{\xyJonesWenzlIdempotentSeven}[1]{\xybox{(0,-9.)*\xybox{
(1.5,0);(22.5,18) **i@{-},(3,18);(3,14.4) **@{-}, (3,3.6);(3,0) **@{-},(6,18);(6,14.4) **@{-}, (6,3.6);(6,0) **@{-},(9,18);(9,14.4) **@{-}, (9,3.6);(9,0) **@{-},(15,18);(15,14.4) **@{-}, (15,3.6);(15,0) **@{-},(18,18);(18,14.4) **@{-}, (18,3.6);(18,0) **@{-},(21,18);(21,14.4) **@{-}, (21,3.6);(21,0) **@{-},(12,16.2) *h+++{\cdots},(12,1.8) *+++{\cdots},(2.25,3.6);(21.75,14.4) **\frm{-}, (12.,9.) *+++{#1}}}}\newcommand{\xyWenzlRecurrenceLastTermSixPlusOne}[1]{\xybox{(0,-20.)*\xybox{*\xybox{*\xybox{
(1.5,0);(19.5,18) **i@{-},(3,18);(3,14.4) **@{-}, (3,3.6);(3,0) **@{-},(6,18);(6,14.4) **@{-}, (6,3.6);(6,0) **@{-},(9,18);(9,14.4) **@{-}, (9,3.6);(9,0) **@{-},(15,18);(15,14.4) **@{-}, (15,3.6);(15,0) **@{-},(18,18);(18,14.4) **@{-}, (18,3.6);(18,0) **@{-},(12,16.2) *h+++{\cdots},(12,1.8) *+++{\cdots},(2.25,3.6);(18.75,14.4) **\frm{-}, (10.5,9.) *+++{#1}}!R!(-1.5,0)*\xybox{
(1.5,0);(4.5,18) **i@{-},(3,18);(3,0)  **\crv{(3,5.4)&(3,12.6)},
}!R!(-1.5,0)}!U*\xybox{*\xybox{
(1.5,0);(10.5,4) **i@{-},(3,4);(3,0)  **\crv{(3,1.2)&(3,2.8)},
(6,4);(6,0)  **\crv{(6,1.2)&(6,2.8)},
(9,4);(9,0)  **\crv{(9,1.2)&(9,2.8)},
}!R!(-1.5,0)*\xybox{*i\xybox{
(1.5,0);(4.5,4) **i@{-},(3,4);(3,0)  **\crv{(3,1.2)&(3,2.8)},
}}!R!(-1.5,0)*\xybox{
(1.5,0);(10.5,4) **i@{-},(3,4);(3,0)  **\crv{(3,1.2)&(3,2.8)},
(6,4);(9,4)  **\crv{(6,2.2)&(9,2.2)},
(9,0);(6,0)  **\crv{(9,1.8)&(6,1.8)},
}!R!(-1.5,0)}!U*\xybox{*\xybox{
(1.5,0);(19.5,18) **i@{-},(3,18);(3,14.4) **@{-}, (3,3.6);(3,0) **@{-},(6,18);(6,14.4) **@{-}, (6,3.6);(6,0) **@{-},(9,18);(9,14.4) **@{-}, (9,3.6);(9,0) **@{-},(15,18);(15,14.4) **@{-}, (15,3.6);(15,0) **@{-},(18,18);(18,14.4) **@{-}, (18,3.6);(18,0) **@{-},(12,16.2) *h+++{\cdots},(12,1.8) *+++{\cdots},(2.25,3.6);(18.75,14.4) **\frm{-}, (10.5,9.) *+++{#1}}!R!(-1.5,0)*\xybox{
(1.5,0);(4.5,18) **i@{-},(3,18);(3,0)  **\crv{(3,5.4)&(3,12.6)},
}!R!(-1.5,0)}!U}}}

\DeclareRobustCommand*{\nicefrac}{\@UnitsNiceFrac}%
\makeatother

\begin{document}

\title{A field guide to categories with $A_n$ fusion rules}
\author{Cain Edie-Michell}
\email{cain.edie-michell@anu.edu.au}
\author{Scott Morrison}
\email{scott.morrison@anu.edu.au}
\address{Mathematical Sciences Institute, Australian National University}
\begin{abstract}
We collate information about the fusion categories with $A_n$ fusion rules. 
This note includes the classification of these categories, a realisation via the Temperley-Lieb categories, 
the auto-equivalence groups (both braided and tensor), identifications of the subcategories of invertible objects, and explicit descriptions of the Drinfeld centres.

The first section describes the classification of these categories (as monoidal, dagger, pivotal, and braided categories). The second section describes the properties of these categories. 
\end{abstract}
\maketitle

\section{Classification}
We say that a fusion category has fusion rules $A_{k+1}$ if there are $k+1$ isomorphism classes of simple objects, which we may index $X_i$ for $i = 0, \ldots, k$, such that:
\begin{itemize}
\item $X_0$ is the tensor identity,
\item $X_i \otimes X_1 \cong X_{i-1} \oplus X_{i+1}$ when $1 \leq i < k$, 
\item $X_k \otimes X_1 \cong X_{k-1}$.
\end{itemize}
The fusion graph for tensoring by $X_{1}$ is thus the $A_{k+1}$ Dynkin diagram:
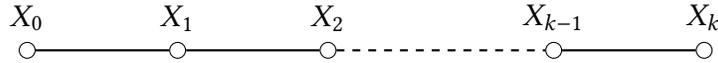
\[
\begin{tikzpicture}
        \vertex[label=$X_0$](f0) at (0,0) {};
        \vertex[label=$X_1$](f1) at (2,0) {};
        \vertex[label=$X_2$](f2) at (4,0) {};
         \vertex[label=$X_{k-1}$](f-2) at (7,0) {};
         \vertex[label=$X_k$](f-1) at (9,0) {};
    %\tikzset{EdgeStyle/.style={-}}
        \Edge(f0)(f1)
        \Edge(f1)(f2)
        \Edge(f-2)(f-1)
     \tikzset{EdgeStyle/.style={dashed}}
        \Edge (f2)(f-2)      
\end{tikzpicture}
\]
 These hypotheses then ensure that the fusion rules are commutative, and in particular that 
\[ X_i \otimes X_j \cong \begin{cases} 
        X_{|i-j|} \oplus X_{|i-j|+2} \oplus \cdots \oplus X_{i+j}   & i+j \leq k\\
         X_{|i-j|} \oplus X_{|i-j|+2} \oplus \cdots \oplus X_{2k -(i+j)} & i+j > k.
   \end{cases}
\]

The only category with $A_1$ fusion rules is $\Vec$. For all $k \geq 1$, the monoidal categories with fusion rules $A_{k+1}$ are classified by $m \in \mathbb{Z} / (k+2)\mathbb{Z}$ with $(m,k+2)=1$. See \cite{MR1239440} for details. We name these monoidal categories $\cC_{k,m}$. The monoidal category $\cC_{k,m}$ has minimal field of definition $\mathbb{Q}(q)$ where $q = e^\frac{m \pi i}{k+2}$. The value $q+q^{-1}$ is a complete invariant of the monoidal categories with $A_{k+1}$ fusion rules.

The monoidal category $\cC_{k,m}$ is generated by an object $X$, and morphisms $\ncup \in \Hom(\mathbf{1} \to X \otimes X)$, $\ncap \in \Hom(X\otimes X \to \mathbf{1})$ satisfying the following relations: 

\begin{align*}
\begin{tikzpicture}
\draw [thick] (0,0) circle (.5);
\end{tikzpicture} & = q+q^{-1}&
\begin{tikzpicture}
\draw [thick] (0,0) -- (0,1);
\draw [thick] (0,0) arc [radius=-.25,start angle = 180, end angle = 0];
\draw [thick] (-.5,0) -- (-.5,1.25);
\draw [thick] (0,1) arc [radius=.25,start angle = 180, end angle = 0];
\draw [thick] (.5,1) -- (.5,-.25);
\end{tikzpicture}=\begin{tikzpicture}
\draw [thick] (0,0) -- (0,1);
\draw [thick] (0,1) arc [radius=.25,start angle = 0, end angle = 180];
\draw [thick] (-.5,1) -- (-.5,-.25);
\draw [thick] (0,0) arc [radius=-.25,start angle = 00, end angle = 180];
\draw [thick] (.5,0) -- (.5,1.25);
\end{tikzpicture}=\begin{tikzpicture}
\draw [thick] (0,-.25) -- (0,1.25);
\end{tikzpicture}&\qquad\qquad\quad f^{(k+1)} = 0.
\end{align*}
Here $f^{(k+1)}$ is the $(k+1)$-th Jones-Wenzl projection, defined recursively so $\JW{0}$ is the identity on $\mathbf{1}$ and
\begin{equation}
\label{eq:wenzl-recursion}
\xyJonesWenzlIdempotentSeven{\JW{k+1}} =
\xyJonesWenzlIdempotentSixPlusOne{\JW{k}} - \frac{[k]_q}{[k+1]_q	}
\xyWenzlRecurrenceLastTermSixPlusOne{\JW{k}}.
\end{equation}

\subsection*{Dagger structures} 
The monoidal category $\cC_{k,m}$ admits a dagger structure $\dagger$ given by the conjugate linear extension of reflection in a horizontal line. All other dagger structures are of the form $\phi \mapsto \lambda^{(n-m)/2} \phi^\dagger$ (where $\phi:\Hom(X^{\otimes n} \to X^{\otimes m})$) for some $\lambda \in \bbR^\times$. We denote this dagger category $\cC_{k,m,\lambda}^\dagger$.

\subsection*{Braidings}
We recall that for $k>1$, the braided fusion categories with fusion rules $A_{k+1}$ are classified by $\ell \in \bbZ / 4(k+2) \bbZ$ with $(\ell, k+2) = 1$. See \cite[Proposition 8.2.6]{MR1239440} for details. We name these braided categories $\cC_{k,l}^\text{br}$. We can recover the underlying monoidal category of $\cC_{k,\ell}^\text{br}$ by taking $m:= \ell\pmod{k+2}$. The braided category $\cC_{k,l}^\text{br}$ has a split rational form over $\bbQ(s)$, where $s = \exp\left(2 \pi i \frac{\ell+k+2}{4(k+2)}\right) = i \exp\left(\pi i \frac{\ell}{2(k+2)}\right)$. We have $s^2 = \pm q$, and this $s$ is a complete invariant of the braided category.

The above classification does not hold for $k=1$ as it double counts the braided fusion categories. The $4$ braided fusion categories with the $A_2$ fusion rules are usually called $\Rep(\mathbb{Z} / 2\mathbb{Z})$, $\sVec$, $\Sem$, and $\overline{\Sem}$. (The category $\Sem$ is the ``semions'', and $\overline{\Sem}$ its Galois conjugate.) The $A_2$ categories can be identified as:
\begin{align*}
\Rep(\mathbb{Z} / 2\mathbb{Z})
	& = \cC_{1,7}^\text{br} =  \cC_{1,11}^\text{br} \\
\sVec				    
	& = \cC_{1,1}^\text{br} =  \cC_{1,5}^\text{br} \\
\Sem      				    
	& = \cC_{1,4}^\text{br} =  \cC_{1,8}^\text{br} \\
\overline{\Sem}                             
	& = \cC_{1,2}^\text{br} =  \cC_{1,10}^\text{br}. 
\end{align*} 

The braided category $\cC_{k,\ell}^\text{br}$ can be realised by equipping the monoidal category $\cC_{k,\ell \pmod {k+2}}$ with a braiding that depends of $\ell$. This braiding is:
\begin{align*}
\begin{tikzpicture}
\draw [thick](0,1) -- (1,0);
\draw [thick,knot,knot gap=7] (0,0) -- (1,1);
\end{tikzpicture} & = s\begin{tikzpicture}
\draw [thick](0,0) arc [radius  = -.707,start angle = 135, end angle = 225];
\draw [thick](1,0) arc [radius  = .707,start angle = 225, end angle = 135];
\end{tikzpicture} + s^{-1}\begin{tikzpicture}
\draw [thick](0,0) arc [radius  = .707,start angle = 135, end angle = 45];
\draw [thick](0,1) arc [radius  = .707,start angle = -135, end angle = -45];
\end{tikzpicture}.
\end{align*}

\subsection*{Monoidal equivalence} 
Recall that the braided category $\cC_{k,\ell}^\text{br}$ has underlying monoidal category $\cC_{k,\ell \pmod{k+2}}$. Thus the braided category $\cC_{k,\ell_1}^\text{br}$ is monoidally equivalent to $\cC_{k,\ell_2}^\text{br}$ if and only if $\ell_1 \equiv \pm \ell_2 \pmod {2(k+2)}$. When $k>1$ each monoidal category admits $4$ distinct braided strucures, and when $k=1$ each monoidal category admits $2$ distinct braided structures.

\subsection*{Pivotality}
The monoidal categories $\cC_{k,m}$ are rigid, with $Y^* = Y$ for every object $Y$, and the coevaluation and evaluation maps are the cups and caps.

The monoidal categories $\cC_{k,m}$ have the rather uninteresting pivotal structure given by the identity morphism on $X$, thought of as an isomorphism $X \iso X^{**} = X$. As $\cC_{k,m}$ has universal grading group $\mathbb{Z} /2\mathbb{Z}$ there exist precisely two pivotal structures \cite{MR2609644}. The other pivotal structure on $\cC_{k,m}$ is given by minus the identity on $X$. We call these two pivotal categories $\cC_{k,m,+}$ and $\cC_{k,m,-}$ respectively. 

There is an alternative skein theoretic presentation of the pivotal categories $\cC_{k,m,-}$, which is generated by two objects $X$ and $X^*$, written as upwards and downwards oriented strands, and morphisms
\begin{align*}
\lcap & \in \Hom( X \otimes X^* \to \mathbf{1}) &
\rcap & \in \Hom( X^* \otimes X \to \mathbf{1})  \displaybreak[1] \\
\lcup & \in \Hom( \mathbf{1} \to X \otimes X^*) &
\rcup & \in \Hom( \mathbf{1} \to X^* \otimes X) \displaybreak[1]  \\
\begin{tikzpicture}
\draw[thick,mid>] (0,0) to (0,.5);
\draw[thick,mid>] (0,1) to (0,.5);
\draw[thick,tag] (0,0.4) -- (0,0.6);
\end{tikzpicture}& \in \Hom( X \to X^*) &
\begin{tikzpicture}
\draw[thick,mid>] (0,0) to (0,.5);
\draw[thick,mid>] (0,1) to (0,.5);
\draw[thick,tag] (0,0.6) -- (0,0.4);
\end{tikzpicture}& \in \Hom( X \to X^*) \displaybreak[1] \\
\begin{tikzpicture}
\draw[thick,mid<] (0,0) to (0,.5);
\draw[thick,mid<] (0,1) to (0,.5);
\draw[thick,tag] (0,0.4) -- (0,0.6);
\end{tikzpicture}& \in \Hom( X^* \to X) &
\begin{tikzpicture}
\draw[thick,mid<] (0,0) to (0,.5);
\draw[thick,mid<] (0,1) to (0,.5);
\draw[thick,tag] (0,0.6) -- (0,0.4);
\end{tikzpicture}& \in \Hom( X^* \to X) 
\end{align*}
satisfying the following relations:
\begin{align*}
\begin{tikzpicture}
\draw [thick,mid>] (0,0) circle (.5);
\end{tikzpicture}  &= 
\begin{tikzpicture}
\draw [thick,mid<] (0,0) circle (.5);
\end{tikzpicture}  = - (q+ q^{-1})
& f^{(k+1)} & = 0  \displaybreak[1] \\
\begin{tikzpicture}
\begin{scope}[very thick, every node/.style={sloped,allow upside down}]
\draw [thick] (0,.25) -- (0,.75);
\draw [thick] (0,.25) arc [radius=-.25,start angle = 180, end angle = 0];
\draw [thick,mid>] (-.5,0.25) -- (-.5,1.25);
\draw [thick] (0,.75) arc [radius=.25,start angle = 180, end angle = 0];
\draw [thick,mid<] (.5,.75) -- (.5,-.25);
\end{scope}
\end{tikzpicture} &=  \begin{tikzpicture}
\draw [thick, mid<] (0,1.25) -- (0,-0.25);
\end{tikzpicture} = \begin{tikzpicture}
\begin{scope}[very thick, every node/.style={sloped,allow upside down}]
\draw [thick] (0,0.25) -- (0,.75);
\draw [thick] (0,.25) arc [radius=-.25,start angle = 0, end angle = 180];
\draw [thick,mid>] (.5,0.25) -- (.5,1.25);
\draw [thick] (0,.75) arc [radius=.25,start angle = 0, end angle = 180];
\draw [thick,mid<] (-.5,.75) -- (-.5,-.25);
\end{scope}
\end{tikzpicture} &
\begin{tikzpicture}
\begin{scope}[very thick, every node/.style={sloped,allow upside down}]
\draw [thick] (0,0.25) -- (0,.75);
\draw [thick] (0,0.25) arc [radius=-.25,start angle = 180, end angle = 0];
\draw [thick,mid<] (-.5,0.25) -- (-.5,1.25);
\draw [thick] (0,.75) arc [radius=.25,start angle = 180, end angle = 0];
\draw [thick,mid>] (.5,.75) -- (.5,-.25);
\end{scope}
\end{tikzpicture} &=  \begin{tikzpicture}
\draw [thick, mid>] (0,1.25) -- (0,-0.25);
\end{tikzpicture} =
\begin{tikzpicture}
\begin{scope}[very thick, every node/.style={sloped,allow upside down}]
\draw [thick] (0,.25) -- (0,.75);
\draw [thick] (0,0.25) arc [radius=-.25,start angle = 0, end angle = 180];
\draw [thick,mid<] (.5,0.25) -- (.5,1.25);
\draw [thick] (0,.75) arc [radius=.25,start angle = 0, end angle = 180];
\draw [thick,mid>] (-.5,.75) -- (-.5,-.25);
\end{scope}
\end{tikzpicture} \\
\begin{tikzpicture}
\begin{scope}[very thick, every node/.style={sloped,allow upside down}]
\draw [thick,tag] (0,0.25) -- (0,.75);
\draw [thick] (0,0.25) arc [radius=-.25,start angle = 180, end angle = 0];
\draw [thick,mid>] (-.5,0.25) -- (-.5,1.25);
\draw [thick] (0,.75) arc [radius=.25,start angle = 180, end angle = 0];
\draw [thick,mid>] (.5,.75) -- (.5,-.25);
\end{scope}
\end{tikzpicture} &=  \begin{tikzpicture}
\draw[thick,mid<] (0,-.25) to (0,.5);
\draw[thick,mid<] (0,1.25) to (0,.5);
\draw[thick,tag] (0,0.6) -- (0,0.4);
\end{tikzpicture}= \begin{tikzpicture}
\begin{scope}[very thick, every node/.style={sloped,allow upside down}]
\draw [thick,tag] (0,0.25) -- (0,.75);
\draw [thick] (0,0.25) arc [radius=-.25,start angle = 0, end angle = 180];
\draw [thick,mid>] (.5,0.25) -- (.5,1.25);
\draw [thick] (0,.75) arc [radius=.25,start angle = 0, end angle = 180];
\draw [thick,mid>] (-.5,.75) -- (-.5,-.25);
\end{scope}
\end{tikzpicture} &
\begin{tikzpicture}
\begin{scope}[very thick, every node/.style={sloped,allow upside down}]
\draw [thick,tag] (0,0.25) -- (0,.75);
\draw [thick] (0,0.25) arc [radius=-.25,start angle = 180, end angle = 0];
\draw [thick,mid<] (-.5,0.25) -- (-.5,1.25);
\draw [thick] (0,.75) arc [radius=.25,start angle = 180, end angle = 0];
\draw [thick,mid<] (.5,.75) -- (.5,-.25);
\end{scope}
\end{tikzpicture} &=  \begin{tikzpicture}
\draw[thick,mid>] (0,-.25) to (0,.5);
\draw[thick,mid>] (0,1.25) to (0,.5);
\draw[thick,tag] (0,0.6) -- (0,0.4);
\end{tikzpicture} = 
\begin{tikzpicture}
\begin{scope}[very thick, every node/.style={sloped,allow upside down}]
\draw [thick,tag] (0,0.25) -- (0,.75);
\draw [thick] (0,0.25) arc [radius=-.25,start angle = 0, end angle = 180];
\draw [thick,mid<] (.5,0.25) -- (.5,1.25);
\draw [thick] (0,.75) arc [radius=.25,start angle = 0, end angle = 180];
\draw [thick,mid<] (-.5,.75) -- (-.5,-.25);
\end{scope}
\end{tikzpicture}  \displaybreak[1] \\
\begin{tikzpicture}
\draw[thick,mid>] (0,0) to (0,.5);
\draw[thick,mid>] (0,1) to (0,.5);
\draw[thick,tag] (0,0.4) -- (0,0.6);
\end{tikzpicture}& = (-1)
\begin{tikzpicture}
\draw[thick,mid>] (0,0) to (0,.5);
\draw[thick,mid>] (0,1) to (0,.5);
\draw[thick,tag] (0,0.6) -- (0,0.4);
\end{tikzpicture} &
\begin{tikzpicture}
\draw[thick,mid<] (0,0) to (0,.5);
\draw[thick,mid<] (0,1) to (0,.5);
\draw[thick,tag] (0,0.4) -- (0,0.6);
\end{tikzpicture} & = (-1)
\begin{tikzpicture}
\draw[thick,mid<] (0,0) to (0,.5);
\draw[thick,mid<] (0,1) to (0,.5);
\draw[thick,tag] (0,0.6) -- (0,0.4);
\end{tikzpicture} \\
\begin{tikzpicture}
\draw[thick,mid<] (0,0) to (0,.5);
\draw[thick,mid<] (0,1) to (0,.5);
\draw[thick,tag] (0,0.6) -- (0,0.4);
\draw[thick,mid<] (0,1) -- (0,1.5);
\draw[thick,tag] (0,.9) -- (0,1.1);
\end{tikzpicture} & = \begin{tikzpicture}
\draw[thick,mid<] (0,0) -- (0,1.5);
\end{tikzpicture} &
\begin{tikzpicture}
\draw[thick,mid>] (0,0) to (0,.5);
\draw[thick,mid>] (0,1) to (0,.5);
\draw[thick,tag] (0,0.6) -- (0,0.4);
\draw[thick,mid>] (0,1) -- (0,1.5);
\draw[thick,tag] (0,.9) -- (0,1.1);
\end{tikzpicture} & = \begin{tikzpicture}
\draw[thick,mid>] (0,0) -- (0,1.5);
\end{tikzpicture}
\end{align*}
Here $f^{(k+1)}$ is defined via a similar recursive formula as in Equation \eqref{eq:wenzl-recursion}: now all strings point upwards, and the two critical points in the second term have tags facing towards each other. In this category $(X^*)^* = X$, and again the pivotal isomorphisms are just the identities. A (non-canonical) monoidal equivalence from $\cC_{k,m,-}$ to $\cC_{k,m,+}$ is given by $X, X^* \mapsto X$ and
\newcommand{\str}{
\begin{tikzpicture}
\draw[thick] (0,0) to (0,1);
\end{tikzpicture}
}
\begin{align}\label{eq:monoidalequivalence}
\lcap & \mapsto \ncap &
\rcap & \mapsto (-1) \;\ncap  &
\lcup & \mapsto \ncup  &
\rcup & \mapsto (-1) \;\ncup  \\
\begin{tikzpicture}
\draw[thick,mid>] (0,0) to (0,.5);
\draw[thick,mid>] (0,1) to (0,.5);
\draw[thick,tag] (0,0.4) -- (0,0.6);
\end{tikzpicture}& \mapsto \str &
\begin{tikzpicture}
\draw[thick,mid>] (0,0) to (0,.5);
\draw[thick,mid>] (0,1) to (0,.5);
\draw[thick,tag] (0,0.6) -- (0,0.4);
\end{tikzpicture}& \mapsto (-1)\;\str  &
\begin{tikzpicture}
\draw[thick,mid<] (0,0) to (0,.5);
\draw[thick,mid<] (0,1) to (0,.5);
\draw[thick,tag] (0,0.4) -- (0,0.6);
\end{tikzpicture}& \mapsto (-1)\;\str  &
\begin{tikzpicture}
\draw[thick,mid<] (0,0) to (0,.5);
\draw[thick,mid<] (0,1) to (0,.5);
\draw[thick,tag] (0,0.6) -- (0,0.4);
\end{tikzpicture}& \mapsto \str \nonumber
\end{align}
See \cite{1002.0555} for a discussion of these two pivotal structures.

The category $\Rep(U_v(\mathfrak{sl}_2))$ has a natural pivotal structure coming from the Hopf algebra structure of the quantum group, for which the dimension of the standard representation is $v+v^{-1}$.
When $v$ is a $2k+4$-th root of unity, the semisimple quotient of $\Rep(U_v(\mathfrak{sl}_2))$ is equivalent as a pivotal category to $\cC_{k,m,-}$ where $v+v^{-1} = -2\cos\left(\frac{m \pi i}{k+2}\right)$. The category of level $k$ integrable highest weight modules of $\widehat{\mathfrak{sl}}_2$ is equivalent as a pivotal category to $\cC_{k,k+1,-}$; this is the category typically referred to as `$\mathfrak{sl}_2$ at level $k$'.

%We write $\cC_{k,\ell,+}^\text{br}$ and $\cC_{k,\ell,-}^\text{br}$ for the corresponding braided pivotal categories. 

\section{Properties}

\subsection*{Basic data}
We now give all the basic data for the braided pivotal categories $\cC_{k,\ell,+}^\text{br}$. The data for $\cC_{k,\ell,-}^\text{br}$ differ by a sprinkling of minus signs as above in Equation~\ref{eq:monoidalequivalence}. Our conventions for quantum integers are $
[n] = (s^{2n} - s^{-2n})/(s^2 - s^{-2})
      %\frac{s^{2n} - s^{-2n}}{s^2 - s^{-2}}
$, and we write $[n]!$ for the quantum factorial $[n][n-1]\cdots [2][1]$. Note that quantum integers only depend on $\delta = q + q^{-1} = -s^2 -s^{-2}$, and hence only on the underlying monoidal category. The categorical dimension of the simple object $X_n$ in the pivotal monoidal category $\cC_{k,m,+}$ is $(-1)^n[n+1]$. The global dimension of $\cC_{k,m,+}$ is thus $\frac{2(k+2)}{2 - s^4 - s^{-4}}$.

Let $a,b,c$ be natural numbers such that the sum of all three is even, and the sum of any two is greater than the third. We say such a triple of numbers is admissible. We choose a basis for the trivalent vertex space of an admissible triple as follows: 

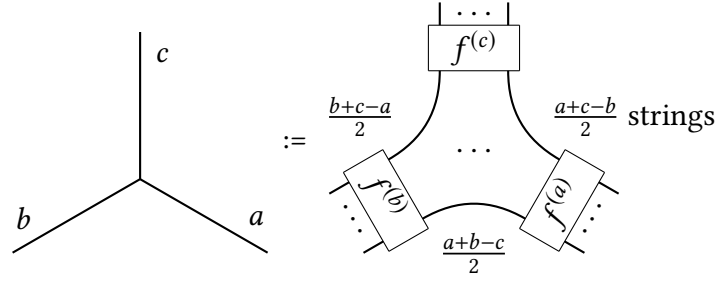
\begin{center}
\begin{tikzpicture}[scale = .303,baseline={([yshift=-.5ex]current bounding box.center)}]
\draw[thick] (0,0) -- (0,6.46);
\node (X) at  (1,5.46) {$c$};
\draw[thick] (0,0) -- (-.866*6.46,-6.46/2);
\node (X) at   (-.866*6.46+.5,-6.46/2+1.5) {$b$};
\draw[thick] (0,0) -- (.866*6.46,-6.46/2);
\node (X) at   (.866*6.46-.5,-6.46/2+1.5) {$a$};
\end{tikzpicture} :=
\begin{tikzpicture}[scale = .303,baseline={([yshift=-.5ex]current bounding box.center)}]
\draw (0,0) -- (4,0);
\draw (0,0) -- (0,-2);
\draw (0,-2) -- (4,-2);
\draw (4,0) -- (4,-2);
\node (X) at  (2,-1) {$f^{(c)}$};

\draw (-2,-3.4641 -2) -- (-2+2,-3.4641 -2  -3.4641);
\draw (-2,-3.4641 -2) -- (-2 - 1.731,-3.4641 -2 - 1);
\draw (-2 - 1.731,-3.4641 -2 - 1) -- (-2 - 1.731+2,-3.4641 -2 - 1-3.4641);
\draw (-2+2,-3.4641 -2  -3.4641)-- (-2 - 1.731+2,-3.4641 -2 - 1-3.4641);
 \node (X) at  (-2,-9) {\raisebox{\dimexpr 2ex}[0pt][0pt]{\rotatebox[origin=c]{-60}{ $f^{(b)}$}}};
 
\draw (2+4,-3.4641 -2) -- (0+4,-3.4641 -2  -3.4641);
\draw (2+4,-3.4641 -2) -- (2 +4+ 1.731,-3.4641 -2 - 1);
\draw (2 +4+ 1.731,-3.4641 -2 - 1) -- (1.731+4,-3.4641 -2 - 1-3.4641);
\draw (0+4,-3.4641 -2  -3.4641)-- (1.731+4,-3.4641 -2 - 1-3.4641);
 \node (X) at  (2+4-.3,-9-.6) {\raisebox{\dimexpr 2ex}[0pt][0pt]{\rotatebox[origin=c]{60}{ $f^{(a)}$}}};
 
 \draw[thick] ( .5,-2) to[bend left = 30] (-2 +.25, -3.4641 -2 - .866*.5);
  \draw[thick] ( 4-.5,-2) to[bend right = 30] (+2 -.25+4, -3.4641 -2 - .866*.5);
\draw[thick]  (-.25, -3.4641 -2 - .866*3.5) to[bend left = 30] (.25+4, -3.4641 -2 - .866*3.5);

\draw[thick] (.5,0) -- (.5,1);
\draw[thick] (3.5,0) -- (3.5,1);

\draw[thick] (-2 - 1.731+.25,-3.4641 -2 - 1- .866*.5) -- (-2 - 1.731+.25-.866,-3.4641 -2 - 1- .866*.5-.5);
\draw[thick] (-2 - 1.731+2 -.25,-3.4641 -2 - 1-3.4641 + .866*.5) -- (-2 - 1.731+2 -.25-.866,-3.4641 -2 - 1-3.4641 + .866*.5-.5);

\draw[thick] (2 + 1.731-.25+4,-3.4641 -2 - 1- .866*.5) -- (2 + 1.731-.25+.866+4,-3.4641 -2 - 1- .866*.5-.5);
\draw[thick] (2 + 1.731-2 +.25+4,-3.4641 -2 - 1-3.4641 + .866*.5) -- (2 + 1.731-2 +.25+.866+4,-3.4641 -2 - 1-3.4641 + .866*.5-.5);

 \node (X) at  (9,-4) {$\frac{a+c-b}{2}$ strings};
  \node (X) at  (-3,-4) {$\frac{b+c-a}{2}$};
   \node (X) at  (2,-10) {$\frac{a+b-c}{2}$};
   
    \node (X) at  (2,-5.5) {$\hdots$};
    
    \node (X) at  (2,.5) {$\hdots$};
     \node (X) at  (-3.4,-9.6) {\raisebox{\dimexpr 2ex}[0pt][0pt]{\rotatebox[origin=c]{-60}{ $\hdots$}}};
     \node (X) at  (4+3.1,-9.9) {\raisebox{\dimexpr 2ex}[0pt][0pt]{\rotatebox[origin=c]{60}{ $\hdots$}}};
\end{tikzpicture}\end{center}
With respect to this choice of basis, the $R$, $S$, and $T$ matrices, and the $\theta$ and $6j$ symbols for $\cC_{k,\ell,+}^\text{br}$ are computed in \cite{MR1280463,MR2640343}.
(Note that the $6j$ symbols only depend on the underlying monoidal category, the $\theta$ symbols depend on the underlying pivotal category, and of course the $R$, $S$, and $T$ matrices depend also on the braiding.)
\begin{align*}
 \theta_{a,b,c} &= \begin{tikzpicture}[scale = .303,baseline={([yshift=-.5ex]current bounding box.center)}]
\draw [thick] (0,0)
	arc[radius = -2, start angle = 0, end angle = 90] 
	node[below] {$c$}
	arc[radius = -2, start angle = 90, end angle = 180] ;
\draw [thick] (0,0)
	arc[radius = -2, start angle = -0, end angle = -90]
	node[above] {$a$}
	arc[radius = -2, start angle = -90, end angle = -180];
\draw [thick] (0,0) -- node[above] {$b$} (4,0);
\end{tikzpicture} =\begin{cases} 
           \frac{ (-1)^{u+v+w}[u+v+w+1]![u]![v]![w]!}{[u+v]![v+w]![u+w]!} &
           \text{if $u,v,$ and $w$ are positive integers}\\
         0 & \text{otherwise}
   \end{cases}\\
 & \qquad \text{where  $u =\frac{b + c-a}{2}$, $v=\frac{ a + c-b}{2}$, $w =\frac{ a + b - c}{2}$}\\
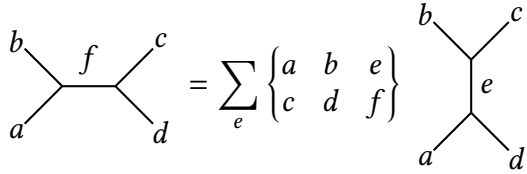
\begin{tikzpicture}[scale = .707,baseline={([yshift=-.5ex]current bounding box.center)}]
\draw [thick] (0,0) -- (-.707,-.707);
\draw [thick] (0,0) -- (-.707,.707);
\draw [thick] (0,0) -- node[above] {$f$} (1,0);
\draw [thick] (1,0) -- (1.707,.707);
\draw [thick] (1,0) -- (1.707,-.707);
\node (X) at  (-.707-.15,-.707-.15) {$a$};
\node (X) at  (-.707-.15,.707+.15) {$b$};
\node (X) at  (1.707+.15,.707+.15) {$c$};
\node (X) at  (1.707+.15,-.707-.15) {$d$};
\end{tikzpicture}&= \sum_{e}\begin{Bmatrix} 
a & b & e \\ 
c & d & f 
\end{Bmatrix} \begin{tikzpicture}[scale = .707,baseline={([yshift=-.5ex]current bounding box.center)}]
\draw [thick] (0,0) -- (-.707,-.707);
\draw [thick] (0,0) -- (.707,-.707);
\draw [thick] (0,0) -- (0,1);
\draw [thick] (0,1) -- (.707,1.707);
\draw [thick] (0,1) -- (-.707,1.707);
\node (X) at  (-.707-.15,-.707-.15) {$a$};
\node (X) at  (.707+.15,-.707-.15) {$d$};
\node (X) at  (.707+.15,1.707+.15) {$c$};
\node (X) at  (-.707-.15,1.707+.15) {$b$};
\node (X) at  (.3,.5) {$e$};
\end{tikzpicture} \\
\begin{Bmatrix} 
a & b & e \\ 
c & d & f
\end{Bmatrix}&=  \begin{cases}
\frac{  \mathcal{I}![e+1](-1)^{e+1}}{\mathcal{E}!\theta_{a,d,e}\theta_{b,c,e}} \sum_{n\leq s \leq N} \frac{ (-1)^s[s+1]!}{\prod_i [s-a_i]!\prod_j[b_j-s]!} & \text{if $(a,d,e)$ is admissible,}\\
0 & \text{otherwise}
\end{cases}
\end{align*}
where
\begin{align*}
\mathcal{I}! &= \prod_{i,j}[b_j-a_i]! & \mathcal{E}! & = [a]![b]![c]![d]![e]![f]! \\
n &= \max\{a_i\} & N & = \min\{b_j\}
\\
% a_1 &= \frac{a+d+i}{2}  &  
% a_2 &= \frac{b+c+i}{2}  &  
% a_3 &= \frac{a+b+j}{2}  &  
% a_4 &= \frac{c+d+j}{2} \\
% b_1 & = \frac{b+d+i+j}{2} &
% b_2 & = \frac{a+c+i+j}{2} &
% b_3 & = \frac{a+b+c+d}{2} &
% \end{align*}
% \begin{align*}
a_1 &= \frac{a+d+e}{2}  & b_1 & = \frac{b+d+e+f}{2} \\
a_2 &= \frac{b+c+e}{2}  & b_2 & = \frac{a+c+e+f}{2} \\
a_3 &= \frac{a+b+f}{2}  & b_3 & = \frac{a+b+c+d}{2} \\
a_4 &= \frac{c+d+f}{2}.
 \end{align*}

\begin{align*}
R_{a,b,c} &=  \begin{tikzpicture}
\draw [thick] (0,0) -- (0,1/2);
\draw [thick] (0,1/2) -- (.293,1.707/2);
\draw [thick] (0,1/2) -- (-.293,1.707/2);
\draw [thick] (-.293,1.707/2) arc [radius = .1,start angle = 225, end angle = 135];
\draw [thick] (.293,1.707/2) arc [radius = .1,start angle = -45, end angle = 45];
\draw [thick] (.293,1.707/2+0.14142135623) -- (-.293,1.707/2+0.14142135623 + .707/2);
\draw [thick,knot,knot gap = 7] (-.293,1.707/2+0.14142135623) -- (.293,1.707/2+0.14142135623 + .707/2);
\node (X) at (.15,0) {$c$};
\node (Y) at (.293,1.707/2+0.14142135623 + .707/2+.2) {$a$};
\node (Z) at (-.293,1.707/2+0.14142135623 + .707/2+.2) {$b$};
\end{tikzpicture} = (-1)^\frac{a+b+c}{2}s^\frac{c(c+2) - a(a+2) - b(b+2)}{2}\begin{tikzpicture}
\draw [thick] (0,0) -- (0,1/2);
\draw [thick] (0,1/2) -- (-.293,1.707/2+0.14142135623 + .707/2);
\draw [thick] (0,1/2) -- (.293,1.707/2+0.14142135623 + .707/2);
\node (X) at (.15,0) {$c$};
\node (Y) at (.293,1.707/2+0.14142135623 + .707/2+.2) {$a$};
\node (Z) at (-.293,1.707/2+0.14142135623 + .707/2+.2) {$b$};
\end{tikzpicture} \\[2ex]
S_{a,b} &=  \begin{tikzpicture}[scale = .5,xscale = -1,baseline={([yshift=-.5ex]current bounding box.center)}]
\draw [thick] (0,0) arc [radius = 1, start angle = 0, end angle = 360];
\draw [thick,draw=white,double=black,double distance=\pgflinewidth,ultra thick] (-.5,.707) arc [radius = 1, start angle = 135, end angle = 135+360 - 15] ;
\node (X) at (.3,0) {$b$};
\node (Y) at (-1.2,0) {$a$};
\end{tikzpicture} =  (-1)^{a+b}[(a+1)(b+1)] \\
T_{a,a} &= \begin{tikzpicture}[scale = .5,yscale = -1,baseline={([yshift=-.5ex]current bounding box.center)}]
\draw [thick](0,0) arc [radius = .5, start angle = 180, end angle = 180+360 - 35];
\draw [thick] (0,0) -- (0,1.5);
\draw [thick] (0,-.35) -- (0,-1.5 );
\node (X) at (.3,1.5) {$a$};
\end{tikzpicture}=   (-1)^{a}s^{a(a+2)} \begin{tikzpicture}[scale = .5,baseline={([yshift=-.5ex]current bounding box.center)}]
\draw [thick] (0,0) -- (0,3);
\node (X) at (.3,0) {$a$};
\end{tikzpicture}
\end{align*}

\subsection*{Galois conjugates} 
The monoidal categories $\cC_{k,m}$ form two Galois equivalence classes when $k$ is odd, depending on the parity of $m$.  When $k$ is even the monoidal categories lie in a single Galois orbit. 
The braided categories $\cC_{k,\ell}^\text{br}$ form three Galois equivalence classes when $k$ is even, depending on whether $\ell \pmod 4$ lies in $\{0\}$, $\{2\}$, or $\{1,3\}$. When $k$ is even the braided categories lie in a single Galois orbit.

% \subsection*{Unitarity}
% The dagger category $\cC_{k,m,\lambda}^\dagger$ is unitary \nn{should we say $C^*$? should we give a reference to this definition?} if and only if $m \equiv  1 \pmod{k+2}$ and $\lambda > 0$ or \nn{...}. We thus see that every Galois equivalence class contains at least one unitary category. \nn{Say the right thing about the relationship with spherical...}

\subsection*{Sphericality}
The pivotal categories $\cC_{k,m,\pm}$ are always spherical.

\subsection*{Modularity}

When $k$ is odd and $2\nmid \ell$ the rank of the $S$-matrix of $\cC_{k,\ell,\pm}^\text{br}$ is $\frac{k+1}{2}$. For all the other cases the $S$-matrix has full rank and hence the categories are modular.

The conductor of $\cC_{k,\ell,\pm}^\text{br}$ (the order of the $T$-matrix ) depends on $k +\ell \pmod 4$. When $k \neq 1$ the conductor is as follows:
\begin{center}\begin{tabular}{c|c}
$k+\ell \pmod 4$ & Conductor \\
\hline
$0$ & $k+2$  \\
$1$ & $4(k+2) $  \\ 
$2$ & $2(k+2)$  \\
$3$ & $4(k+2)$  
\end{tabular}\end{center}
In particular when $k$ is even the conductor is always $4(k+2)$.

The conductors of $\Rep(\mathbb{Z} / 2\mathbb{Z})$, $\sVec$, $\Sem$, and $\overline{\Sem}$ are $1$,$2$,$4$, and $4$ respectively.

\subsection*{Invertible objects}
There are always two invertible objects in a monoidal category with $A_{k+1}$ fusion rules.
We now describe the rule identifying, as a braided category, the maximal pointed subcategory of the braided category $\cC_{k,l}^\text{br}$. This subcategory must be one of $\Rep(\mathbb{Z} / 2\mathbb{Z}), \sVec, \Sem,$ or $\overline{\Sem}$. (Of course, the monoidal structure on the maximal pointed subcategory only depends on the monoidal structure of $\cC_{k,l}^\text{br}$.)

We can determine which appears by choosing a pivotal structure for $\cC_{k,\ell}^\text{br}$ and considering the quantum twist of the non-trivial object of the $A_2$ subcategory, and comparing to the quantum twists of the $A_2$ categories with the appropriately chosen pivotal structure. We get the following identification:
 
\begin{center}\begin{tabular}{cc|cccc}
$\Inv(\cC_{k,\ell}^\text{br})$ & & & \multicolumn{2}{c}{$k \pmod 4$} \\
& & $0$& $1$& $2$ & $3$  \\ \hline 
&  $0$ & --- &  $\Sem$& --- & $\Sem$  \\[2ex]
&  $1$ &$\Rep(\mathbb{Z} / 2\mathbb{Z})$  &$\sVec$& $\sVec$ &$\Rep(\mathbb{Z} / 2\mathbb{Z})$\\[2ex]
\raisebox{\dimexpr 2ex}[0pt][0pt]{\rotatebox[origin=c]{90}{$\ell \pmod 4$}} &  $2$ & --- &  $\overline{\Sem}$ & --- & $\overline{\Sem}$ \\[2ex]
&  $3$ & $\Rep(\mathbb{Z} / 2\mathbb{Z})$  &$\Rep(\mathbb{Z} / 2\mathbb{Z})$ & $\sVec$& $\sVec$. 
\end{tabular}\end{center}
In particular when $k  \equiv 0 \pmod 4$ the invertible elements of $\cC_{k,\ell}^\text{br}$ always form a $\Rep(\mathbb{Z} / 2\mathbb{Z})$ subcategory, and when $k  \equiv 2 \pmod 4$ the invertible elements of $\cC_{k,\ell}^\text{br}$ always form an $\sVec$ subcategory.

\subsection*{Auto-equivalences}
When $k \geq 3$ the groups of tensor auto-equivalences of $\cC_{k,m}$, and braided auto-equivalences of $\cC_{k,\ell}^\text{br}$ are computed in \cite{1709.04721}. This result is independent of the choice of $m$ and $\ell$. 
\begin{center}\begin{tabular}{c|c|c}
$k\pmod 4$ & $\operatorname{Aut}_\otimes ( \cC_{k,m})$ & $\operatorname{Aut}_\otimes^\text{br} ( \cC_{k,\ell}^\text{br})$ \\
\hline
$0$ & $\mathbb{Z}/2\mathbb{Z}$ & $\{e\}$ \\
$1$ & $\{e\} $ & $\{e\}$ \\ 
$2$ & $\mathbb{Z}/2\mathbb{Z}$ & $\mathbb{Z}/2\mathbb{Z}$ \\
$3$ & $\{e\} $ & $\{e\}$ 
\end{tabular}\end{center}
When $k \leq 2$ both auto-equivalence groups are trivial.

\subsection*{Equivariantisations}
When $k$ is even and at least 4 the monoidal categories $\cC_{k,m}$ can be equivariantised by the $\mathbb{Z} / 2\mathbb{Z}$ action. We can compute the fusion rules for the equivariantisations from the results of \cite{MR3059899}. We give the fusion graph for the generating object $(X_1 \oplus X_{k-1}, \begin{psmallmatrix}0&\operatorname{id}_{X_{k-1}}\\\operatorname{id}_{X_{1}}&0\end{psmallmatrix})$. 

When $k \equiv  0\pmod 4$ the fusion graph is
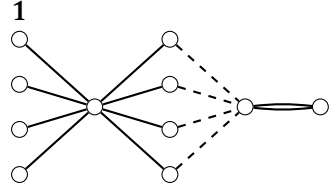
\[
\begin{tikzpicture}[yscale=0.6]
        \vertex[label=$\mathbf{1}$](g1) at (0,3) {};
        \vertex[](g2) at (0,2) {};
        \vertex[](g3) at (0,1) {};
         \vertex[](g4) at (0,0) {};
         
         \vertex[](X) at (1,1.5) {};
         
         \vertex[](Y1) at (2,3) {};
        \vertex[](Y2) at (2,2) {};
        \vertex[](Y3) at (2,1) {};
         \vertex[](Y4) at (2,0) {};
         
          \vertex[](Z) at (3,1.5) {};
          
          \vertex[](W) at (4,1.5) {};
    %\tikzset{EdgeStyle/.style={-}}
        \Edge(g1)(X)
        \Edge(g2)(X)
        \Edge(g3)(X)
         \Edge(g4)(X)
         
          \Edge(Y1)(X)
        \Edge(Y2)(X)
        \Edge(Y3)(X)
         \Edge(Y4)(X)
         
         \Edge[style = {bend left= 10}](W)(Z)
         \Edge[style = {bend right= 10}](W)(Z)
     \tikzset{EdgeStyle/.style={dashed}}
         \Edge(Y1)(Z)
        \Edge(Y2)(Z)
        \Edge(Y3)(Z)
         \Edge(Y4)(Z)    
\end{tikzpicture}
\]

When $k \equiv  2\pmod 4$ the fusion graph is 
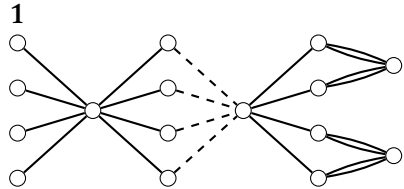
\[\quad\quad
\begin{tikzpicture}[yscale=0.6]
        \vertex[label=$\mathbf{1}$](g1) at (0,3) {};
        \vertex[](g2) at (0,2) {};
        \vertex[](g3) at (0,1) {};
         \vertex[](g4) at (0,0) {};
         
         \vertex[](X) at (1,1.5) {};
         
         \vertex[](Y1) at (2,3) {};
        \vertex[](Y2) at (2,2) {};
        \vertex[](Y3) at (2,1) {};
         \vertex[](Y4) at (2,0) {};
         
          \vertex[](Z) at (3,1.5) {};
          
           \vertex[](W1) at (4,3) {};
        \vertex[](W2) at (4,2) {};
        \vertex[](W3) at (4,1) {};
         \vertex[](W4) at (4,0) {};
         
         \vertex[](U1) at (5,2.5) {};
         \vertex[](U2) at (5,.5) {};
    %\tikzset{EdgeStyle/.style={-}}
        \Edge(g1)(X)
        \Edge(g2)(X)
        \Edge(g3)(X)
         \Edge(g4)(X)
         
          \Edge(Y1)(X)
        \Edge(Y2)(X)
        \Edge(Y3)(X)
         \Edge(Y4)(X)
         
          \Edge(W1)(Z)
        \Edge(W2)(Z)
        \Edge(W3)(Z)
         \Edge(W4)(Z)
         
         \Edge[style = {bend right= 10}](W1)(U1)
         \Edge[style = {bend left= 10}](W1)(U1)
         \Edge[style = {bend right= 10}](W2)(U1)
         \Edge[style = {bend left= 10}](W2)(U1)
         
         \Edge[style = {bend right= 10}](W3)(U2)
         \Edge[style = {bend left= 10}](W3)(U2)
         \Edge[style = {bend right=10}](W4)(U2)
         \Edge[style = {bend left = 10}](W4)(U2)
     \tikzset{EdgeStyle/.style={dashed}}
         \Edge(Y1)(Z)
        \Edge(Y2)(Z)
        \Edge(Y3)(Z)
         \Edge(Y4)(Z)    
\end{tikzpicture}
\]  The depth of the graph in either case is $\frac{k}{2}+1$.

\subsection*{Drinfeld centre}
The monoidal category $\cC_{k,m}$ can be equiped with a braiding to give the braided category $\cC_{k,m}^\text{br}$. This braided category is modular when $k$ is even or $2 \mid m$. For these cases we can use the results of \cite{MR1990929} to see that the Drinfeld centre of the monoidal category $\cC_{k,m}$ is the braided category $\cC_{k,m}^\text{br} \boxtimes \cC_{k,m}^\text{br op}$.

When $k$ is odd and $2\nmid m$ the monoidal category $C_{k,m}$ does not admit a modular braiding, so we have to work harder to compute the centre. As monoidal categories we have an equivalence $C_{k,m} \cong \ad(C_{k,m}) \boxtimes \cC_{1,1}$. Thus $Z(C_{k,m}) \cong Z(\ad(C_{k,m}))\boxtimes Z(\cC_{1,1})$ as the Drinfeld centre respects products of categories \cite{1507.00503}.

The monoidal category $\ad(C_{k,m})$ admits a modular braiding, that is the braiding inherited from the braided category $\cC_{k,m}^\text{br}$. Thus $Z(\ad(C_{k,m})) \cong \ad(C_{k,m}^\text{br})\boxtimes \ad(C_{k,m}^\text{br op})$.

The monoidal category $\cC_{1,1}$ is better known as $\Vec(\mathbb{Z}/2\mathbb{Z})$. The Drinfeld centre of this monoidal category is described in \cite{MR1976233} and a skein theoretic presentation is given in \cite{1709.04721}.

\subsection*{Algebra objects}
The simple algebra objects $A$ in the monoidal categories $\cC_{k,m}$ are as follows.

\begin{center}\begin{tabular}{c|ccc}
$A$ &   when &  $A$-mod  & \text{commutative in $\cC_{k,\ell}^\text{br}$}  \\[.1cm] \hline 
$\mathbf{1} \oplus f^{(k)}$ & $k$ even & $D_{\frac{k}{2}+2}$ & 
				if $k \equiv 0 \pmod 4$ 
			   \\[.5cm]
$\mathbf{1} \oplus f^{(k)}$ & $k,m$ both odd &  $ T_\frac{k+1}{2}$ & 
			  if $k\ell \equiv 3 \pmod 4$ 
			  \\[.5cm]
$\mathbf{1} \oplus f^{(6)}$ & $k = 10$ & $E_6$& yes \\[.5cm]
$\mathbf{1} \oplus f^{(8)} \oplus f^{(16)}$ & $k = 16$ & $E_7$&  no    \\[.5cm]
$\mathbf{1} \oplus f^{(10)}\oplus f^{(18)} \oplus f^{(28)}$ & $k = 28$&  $E_8$ &yes
\end{tabular}\end{center}
See \cite{MR1976459} for details. For details on the categories $A$-mod (in subfactor language) see \cite{MR1193933,MR1145672,MR1313457,MR1929335,MR1308617, MR1617550}. When $A$ is commutative the category $A$-mod has the structure of a fusion category. For details on the skein theory of these fusion categories see \cite{MR2577673,MR2559686}

\subsection*{Summary}

We conclude this survey with a table summarising some of the data of the braided categories $\cC_{k,\ell}^\text{br}$ described earlier.
\vspace{1cm}
\bgroup
\def\arraystretch{2.25}
{\begin{center}\begin{tabular}{ccccc}
& Galois orbits &  modular &  & $\Inv(\cC_{k,\ell}^\text{br})$ \\
\hline\hline
$k$ even
		& a single orbit 			  
			& yes 
			& 
			& \parbox{6cm}{
				$\Rep(\mathbb{Z}/2\mathbb{Z})$ if $k \equiv 0 \pmod 4$ \\
				$\sVec$ \quad \quad \quad if $k \equiv 2 \pmod 4$
			  }
			\\ \cline{2-5}
\multirow{3}{*}{$k$ odd}
		& $\ell \equiv 0 \pmod 2$ 
			& yes
			&  
			& \parbox{6cm}{
                                           $\Sem$ if $\ell \equiv 0 \pmod 4$   \\  
 				$\overline{\Sem}$ if $\ell \equiv 2 \pmod 4$
			}
			\\ \cline{3-5}
		& $\ell \equiv 1 \pmod 4$ 
			& no
			&  
			& 
			\parbox{6cm}{
				$\Rep(\mathbb{Z}/2\mathbb{Z})$ if $k \equiv 3 \pmod 4$ \\
				$\sVec$ \quad \quad \quad if $k \equiv 1 \pmod 4$
			  }
			\\ \cline{3-5}
		& $\ell \equiv 3 \pmod 4$	
			& no
			& 
			&\parbox{6cm}{
				$\Rep(\mathbb{Z}/2\mathbb{Z})$ if $k \equiv 1 \pmod 4$ \\
				$\sVec$ \quad \quad \quad if $k \equiv 3 \pmod 4$
			  }			
			\end{tabular}\end{center}}

\vspace{1cm}
\egroup

%%%%%%%%%%%%%%%%%%%%%%%%%%%%%%%%%%%%%%%%%%%%%%%%%%%%%%%%%%%%
\subsection*{Acknowledgements}
The authors were partially supported by Australian Research Council grants
`Subfactors and symmetries' DP140100732 and `Low dimensional categories' DP160103479. The first author was supported by an Australian Government Research Training Program (RTP) Scholarship.

\renewcommand*{\bibfont}{\small}
\setlength{\bibitemsep}{0pt}
\raggedright

\printbibliography
\end{document}